\documentclass[12pt,reqno]{article}

\usepackage[usenames]{color}
\usepackage{amssymb}
\usepackage{amsmath}
\usepackage{amsthm}
\usepackage{amsfonts}
\usepackage{amscd}
\usepackage{graphicx}

\usepackage[colorlinks=true,
linkcolor=webgreen,
filecolor=webbrown,
citecolor=webgreen]{hyperref}

\definecolor{webgreen}{rgb}{0,.5,0}
\definecolor{webbrown}{rgb}{.6,0,0}

\usepackage{color}
\usepackage{fullpage}
\usepackage{float}

\usepackage{graphics}
\usepackage{latexsym}
\usepackage{epsf}
\usepackage{breakurl}

\def\modd#1 #2{#1\ \mbox{\rm (mod}\ #2\mbox{\rm )}}
\usepackage[shortlabels]{enumitem}

\begin{document}

\theoremstyle{plain}
\newtheorem{theorem}{Theorem}
\newtheorem{corollary}[theorem]{Corollary}
\newtheorem{lemma}[theorem]{Lemma}
\newtheorem{proposition}[theorem]{Proposition}

\theoremstyle{definition}
\newtheorem{definition}[theorem]{Definition}
\newtheorem{example}[theorem]{Example}
\newtheorem{conjecture}[theorem]{Conjecture}

\theoremstyle{remark}
\newtheorem{remark}[theorem]{Remark}

\begin{center}
\vskip 1cm{\LARGE\bf 
Generalization of the Extended Minimal \\
\vskip .1in
Excludant of Andrews and Newman
}
\vskip 1cm
\large
Aritram Dhar, Avi Mukhopadhyay, and Rishabh Sarma\\
Department of Mathematics\\
University of Florida\\
Gainesville, FL 32611\\
USA\\
\href{mailto:aritramdhar@ufl.edu}{\tt aritramdhar@ufl.edu} \\
\href{mailto:mukhopadhyay.avi@ufl.edu}{\tt mukhopadhyay.avi@ufl.edu} \\
\href{mailto:rishabh.sarma@ufl.edu}{\tt rishabh.sarma@ufl.edu} 
\end{center}

\vskip .2in

\begin{abstract}
In a recent pioneering work, Andrews and Newman defined an extended function $p_{A,a}(n)$ of their minimal excludant or ``mex" of a partition function. 
By considering the special cases $p_{k,k}(n)$ and $p_{2k,k}(n)$, they unearthed connections to the rank and crank of partitions and some restricted partitions. In this paper, we build on their work and obtain more general results associating the extended mex function with the number of partitions of an integer with arbitrary bound on the rank and crank. We also derive a new result expressing the smallest parts function of Andrews as a finite sum of the extended mex function in consideration with a curious coefficient. We also obtain a few restricted partition identities with some reminiscent of shifted partition identities. Finally, we define and explore a new minimal excludant for overpartitions.
\end{abstract}

\section{Introduction}\label{intro}
In a series of two papers \cite{An-New19, An-New20}, Andrews and Newman introduced the minimal excludant or mex of an integer partition function and initiated the study of the connections of the mex and its extended functions to other partition statistics. Given a partition $\pi$ of an integer $n$, they defined mex($\pi$) to be the smallest positive integer that is not a part of $\pi$ \cite{An-New19}.

Let $p(n)$ denote the number of partitions of $n$. Following Andrews and Newman \cite{An-New20}, we define mex$_{A,a}(\pi)$ to be the smallest integer greater than or equal to $a$ and  congruent to $a$ modulo $A$ that is not a part of the partition $\pi$. This restriction \cite{An-New20} on the smallest integer being greater than or equal to $a$ is not explicitly stated, but we believe is implied. Further, define $p_{A,a}(n)$ to be the number of partitions $\pi$ of $n$, where $\text{mex}_{A,a}(\pi) \equiv \modd{a} {2A}$ and $\overline{p}_{A,a}(n)$ to be the number of partitions $\pi$ of $n$, where $\text{mex}_{A,a}(\pi) \equiv \modd{A+a} {2A}$. Then $p(n)=p_{A,a}(n)+ \overline{p}_{A,a}(n)$. In Table \ref{p236}, we illustrate an example when $a>A$ for the partitions of $6$. 

\begin{table}[htb]	
\centering
\scalebox{0.9}{%
\begin{tabular}{|l|c|c|c|c|c|}

\hline
Partition $\pi$ & mex$_{2,3}(\pi)$ & $p_{2,3}(6)$ & $\overline{p}_{2,3}(6)$ \\
 \hline
6               & 3               &\checkmark        &      \\
 \hline
5+1             & 3               &\checkmark        &      \\
 \hline
4+2             & 3               &\checkmark        &       \\
 \hline
4+1+1           & 3               &\checkmark        &      \\
 \hline
3+3             & 5               &                  &\checkmark      \\
 \hline
3+2+1           & 5               &                  &\checkmark     \\
 \hline
3+1+1+1         & 5               &                  &\checkmark      \\
 \hline
2+2+2           & 3               &\checkmark        &       \\
 \hline
2+2+1+1         & 3               &\checkmark        &       \\
 \hline
2+1+1+1+1       & 3               &\checkmark        &      \\
 \hline
1+1+1+1+1+1     & 3               &\checkmark        &       \\
 \hline

\end{tabular}}

\caption{Calculating $p_{2,3}(6)$ and $\overline{p}_{2,3}(6)$ from mex$_{2,3}(\pi)$ for partitions $\pi$ of $6$.}\label{p236}

\end{table}

By convention, we define $p_{A,a}(0)=p(0)=1$ and $p(n)=p_{A,a}(n)=\overline{p}_{A,a}(n)=0$ for negative $n$. We also define $$F_{A,a}=\sum\limits_{n\geq 0}p_{A,a}(n)q^n$$ and $$\overline{F}_{A,a}=\sum\limits_{n\geq 0}\overline{p}_{A,a}(n)q^n.$$ 

We use the standard $q$-series notation $(a;q)_n=(1-a)(1-aq)\cdots(1-aq^{n-1})$ and $(a;q)_{\infty}=\lim\limits_{n\rightarrow\infty}(a;q)_n$. When $a=q$, we write $(q;q)_n=(q)_n$ and $(q;q)_{\infty}=(q)_{\infty}$. We also define $p_e(n)$ and $p_o(n)$ to be the number of partitions of $n$ into an even and odd number of parts respectively. Then Andrews and Newman \cite{An-New20} noted that $$\sum\limits_{n\geq 0}p_e(n)q^n=\sum\limits_{n\geq 0}\frac{q^{2n}}{(q)_{2n}}$$ and $$\sum\limits_{n\geq 0}p_o(n)q^n=\sum\limits_{n\geq 0}\frac{q^{2n+1}}{(q)_{2n+1}}.$$

We recall the partition statistics rank, crank, and \text{spt}. The 
{\it rank\/} of a partition is defined as the largest part minus the number of parts. Let $N(m,n)$ be the number of partitions of $n$ with rank $m$. The 
{\it crank\/} of a partition is defined as the largest part of the partition if there are no ones as parts, and otherwise is the number of parts larger than the number of ones minus the numbers of ones. Let $M(m,n)$ be the number of partitions of $n$ with crank $m$. Finally, \text{spt}($n$) is the total number of appearances of the smallest parts in all the partitions of $n$.

Further, the $k$-th Atkin-Garvan rank moment \cite{Atk-Ga03} is defined by $N_k(n)=\sum\limits_{m=-\infty}^{\infty} m^k N(m,n)$. Likewise the $k$-th Atkin-Garvan crank moment \cite{Atk-Ga03} is defined by $M_k(n)=\sum\limits_{m=-\infty}^{\infty} m^k M(m,n)$.

Further, the Garden of Eden partitions of $n$ \cite{Hop-Sel07} are defined as the partitions of $n$ that have rank $-2$ or less.

Andrews and Newman \cite{An-New20} derived the generating functions for $p_{k,k}(n)$ and $p_{2k,k}(n)$ and concluded some striking results, which are as follows:

\begin{theorem}[\cite{An-New20}, Theorem $2$]\label{AnNew1}
If $n$ is a non-negative integer, then $p_{1,1}(n)$ equals the number of partitions of $n$ with non-negative crank.
\end{theorem}

\begin{theorem}[\cite{An-New20}, Theorem $3$]\label{AnNew2}
If $n$ is a non-negative integer, then $p_{3,3}(n)$ equals the number of partitions of $n$ with rank $\geq -1$.
\end{theorem}

\begin{theorem}[\cite{An-New20}, Theorem $4$]\label{AnNew3}
If $n$ is a non-negative integer, then $p_{2,1}(n)=p_e(n)$.
\end{theorem}

In this work, we find the generating function for $p_{A,a}(n)$ in general for arbitrary positive integers $A,a$, which leads us to discover more general results associating the mex with the rank, crank, and spt statistics via elementary techniques, among other restricted partition identities and some auxiliary results. As a consequence, we obtain the above three results of Andrews and Newman as corollaries. Our motivation lies in the questions posed by the authors \cite{An-New20} where they call for the search of comparable theorems and related statistics for other $p_{A,a}(n)$.

We now present the statements of our main results. 

\begin{theorem}\label{mex3132}
For every integer $n\geq 1$, we have $$p_{3,1}(n)+p_{3,2}(n)=p(n).$$
\end{theorem}

\begin{theorem}\label{mexrecur}
If $n$ is a non-negative integer, then $$p_{A,a}(n)=p(n)+\sum\limits_{m\geq 1}(p(n-(A{2m \choose 2}+a(2m)))-p(n-(A{2m-1 \choose 2}+a(2m-1)))).$$
\end{theorem}

\begin{theorem}\label{mexrank1}
If $n$ is a non-negative integer, then $\overline{p}_{3,j+1}(n)$ equals the number of partitions of $n$ with rank $\geq j$, where j is a non-negative integer. 
\end{theorem}

\begin{theorem}\label{mexcrank1}
If $n$ is a non-negative integer, then $\overline{p}_{1,j}(n)$ equals the number of partitions of $n$ with crank $\geq j$, where j is a non-negative integer. 
\end{theorem}

As a consequence of Theorems \ref{mexrank1} and \ref{mexcrank1}, we obtain Corollaries \ref{mexrank2} and \ref{mexcrank2} later, that generalize Theorems \ref{AnNew1} and \ref{AnNew2} of Andrews and Newman, the proofs of which are included in the next section.

\begin{theorem}\label{mexrankmom} 
Let $n$ be a non-negative integer. Then
\begin{align*}
N_2(n)=2\sum\limits_{r=0}^{n-2}(2r+1)\overline{p}_{3,r+2}(n),
\\
M_2(n)=2\sum\limits_{r=0}^{n-1}(2r+1)\overline{p}_{1,r+1}(n).
\end{align*}
\end{theorem}

\begin{theorem}\label{mexspt}
Let $n$ be a non-negative integer. Then
$$\textnormal{spt} (n)=\sum\limits_{r=0}^{n-1}(2r+1)(\overline{p}_{1,r+1}(n)-\overline{p}_{3,r+2}(n))
=\sum\limits_{r=0}^{n-1}(2r+1)(p_{3,r+2}(n)-p_{1,r+1}(n)).$$
\end{theorem}

\begin{theorem}\label{mexcongJTP} (Congruences for odd and even moduli) Let $n$ be a non-negative integer. Then, for every integer $i\geq 1$, we have 
the following results:
\begin{enumerate}[(a)]
\item $p_{2k,2k-i}(n)-\overline{p}_{2k,i}(n)$ equals the number of partitions of $n$ into parts $\not \equiv \modd{0,\pm i} {2k}$ for every integer $k\geq 2$.
\item $p_{2k+1,2k+1-i}(n)-\overline{p}_{2k+1,i}(n)$ equals the number of partitions of $n$ into parts $\not \equiv \modd{0,\pm i} {\\2k+1}$ for every integer $k\geq 1$.
\end{enumerate}
\end{theorem}

As a consequence of the above theorem, we obtain the aforementioned result Theorem \ref{AnNew3} of Andrews and Newman as a corollary, the proof of which is included in the next section.

\begin{theorem}\label{mexcong1}
If $n$ is a non-negative integer, then $p_{2,3}(n)-p_{4,6}(n-1)$ equals the number of partitions of $n$ into parts $\equiv \modd{\pm 2, \pm 8, \pm 12, \pm 14} {32}.$
\end{theorem}

\begin{theorem}\label{mexcong2}
If $n$ is a non-negative integer, then $p_{2,3}(n)-p_{6,9}(n-2)$ equals the number of partitions of $n$ into parts $\equiv \modd{\pm 1, \pm 4, \pm 6, \pm 8, \pm 10, \pm 14} {24}.$
\end{theorem}

\begin{theorem}\label{mexcong3}
If $n$ is a non-negative integer, then $p_{2,3}(n)-p_{10,15}(n-4)$ equals the number of partitions of $n$ into distinct parts $\equiv \modd{\pm8, \pm12} {40}$ and into parts $\not \equiv \modd{0, \pm3, \pm4, \pm7, \pm10,\\ \pm13, \pm17, 20} {40}.$
\end{theorem}
Tables \ref{rc} and \ref{spt} illustrate our results for the rank, crank, and spt for the first few values of $n$ respectively.
\begin{table}[htb]	
\centering
\scalebox{0.85}{
\begin{tabular}{|c|c|c|c|c|c|c|c|c|c|c|c|}
\hline
$n$ & $p_{3,1}(n)$ & $p_{3,2}(n)$   & $\overline{p}_{3,3}(n)$  & $\overline{p}_{1,2}(n)$ & Partitions of $n$  & Partitions of $n$ \\ &&&&& with rank $\geq j=2$ &with crank $\geq j=2$   \\
 \hline
1                  &0                &1               &0              &0               &$\phi$             &$\phi$ \\
 \hline
2                  &1                &1               &0              &1               &$\phi$           &2\\
 \hline
3                  &1                &2               &1              &1               &3             &3\\
 \hline
4                  &2                &3               &1              &2               &4             &4, 2+2\\
 \hline
5                  &3                &4               &2              &2               &5, 4+1             &5, 3+2\\
 \hline
6                  &5                &6               &3              &4               &6, 5+1, 4+2             &6, 4+2, 3+3, 2+2+2\\
 \hline
7                  &6                &9               &5              &5               &7, 6+1, 5+2, 5+1+1, 4+3             &7, 5+2, 4+3, 3+2+2, 2+2+2+1\\
 \hline
8                  &10               &12              &7              &8               &8, 7+1, 6+2, 6+1+1, &8, 6+2, 5+3, 4+4, 4+2+2, \\ &&&&&5+3, 5+2+1, 4+4  &3+3+2, 3+2+2+1, 2+2+2+2\\
 \hline
\end{tabular}}
\caption{Example illustrating Theorems \ref{mex3132}, \ref{mexrank1}, and \ref{mexcrank1}.}\label{rc}
\end{table}

\begin{table}[htb]	
\centering
\scalebox{0.88}{
\begin{tabular}{|c|c|c|c|c|c|c|c|c|c|c|c|}
\hline
$n$ &  $p_{3,2}(n)$ & $p_{1,1}(n)$ & $p_{3,3}(n)$ & $p_{1,2}(n)$ & $p_{3,4}(n)$ & $p_{1,3}(n)$ & $p_{3,5}(n)$ & $p_{1,4}(n)$ & $p_{3,6}(n)$ & $p_{1,5}(n)$ & \text{spt}($n$)\\
 \hline
1                  &1               &0               &1              &1               &1              &1                 &1               &1               &1              &1               &1 \\
 \hline
2                  &1               &1               &2              &1               &2              &2                 &2               &2               &2              &2               &3 \\
 \hline
3                  &2               &2               &2              &2               &3              &2                 &3               &3               &3              &3               &5 \\
 \hline
4                  &3               &3               &4              &3               &4              &4                 &5               &4               &5              &5               &10 \\
 \hline
5                  &4               &4               &5              &5               &6              &5                 &6               &6               &7              &6               &14 \\
 \hline
\end{tabular}}
\caption{Example illustrating Theorem \ref{mexspt}.}\label{spt}
\end{table}
The paper is organized as follows. In Section \ref{proofs}, we give the proofs of our main results. Section \ref{auxiliary} is devoted to studying some additional properties of the extended mex function $p_{A,a}(n)$. We conclude with a glimpse into extending the minimal excludant to overpartitions in Section \ref{overmex}. 
\section{Proof of the main results}\label{proofs}
We simultaneously state and prove some background theorems and lemmas in the build up to the proof of our main results stated in the previous section. In the following two theorems, we find the generating functions for $p_{A,a}(n)$ and $\overline{p}_{A,a}(n)$ for arbitrary values of $A$ and $a$. These two are at the heart of proving all our subsequent results.
\begin{theorem}\label{mexgen1}
$F_{A,a}(q)=\frac{1}{(q)_{\infty}}\sum\limits_{n\geq 0}(-1)^nq^{\frac{An(n-1)}{2}+an}$.
\end{theorem}
\begin{proof}
The generating function for $p_{A,a}(n)$ is 
\begin{align*}
\sum\limits_{n\geq 0}\dfrac{q^{a+(a+A)+\cdots+(a+(2n-1)A)}}{\prod\limits_{\substack{m=1 \\ m\neq a+2nA}}^{\infty}(1-q^m)}&=\frac{1}{(q)_{\infty}}\sum\limits_{n\geq 0}q^{a+(a+A)+\cdots+a+(2n-1)A}(1-q^{a+2nA})
\\
&=\frac{1}{(q)_{\infty}}\sum\limits_{n\geq 0}q^{2n^2A-nA+2na}(1-q^{a+2nA})
\\
&=\frac{1}{(q)_{\infty}}\sum\limits_{n\geq 0}(q^\frac{(2n)^2A-(2n)A+2(2n)a}{2}-q^\frac{(2n+1)^2A-(2n+1)A+2(2n+1)a}{2})
\end{align*}
\begin{align*}
&=\frac{1}{(q)_{\infty}}\sum\limits_{n\geq 0}(-1)^nq^{\frac{n^2A+(2a-A)n}{2}}
\\
&=\frac{1}{(q)_{\infty}}\sum\limits_{n\geq 0}(-1)^nq^{\frac{An(n-1)}{2}+an}.
\end{align*}
\end{proof}
\begin{theorem}\label{mexgen2}
$\overline{F}_{A,a}(q)=\frac{1}{(q)_{\infty}}\sum\limits_{n\geq 0}(-1)^nq^{\frac{An(n+1)}{2}+a(n+1)}$.
\end{theorem}
\begin{proof}
The generating function for $\overline{p}_{A,a}(n)$ is 
\begin{align*}
\sum\limits_{n\geq 0}\dfrac{q^{a+(a+A)+(a+2A)+\cdots+(a+2nA)}}{\prod\limits_{\substack{m=1 \\ m\neq a+(2n+1)A}}^{\infty}(1-q^m)}&=\frac{1}{(q)_{\infty}}\sum\limits_{n\geq 0}q^{a+(a+A)+(a+2A)+\cdots+(a+2nA)}(1-q^{a+(2n+1)A})
\\
&=\frac{1}{(q)_{\infty}}\sum\limits_{n\geq 0}q^{(2n+1)(nA+a)}(1-q^{a+A+2nA})
\\
&=\frac{1}{(q)_{\infty}}\sum\limits_{n\geq 0}(-1)^nq^{\frac{(n+1)(nA+2a)}{2}}
\\
&=\frac{1}{(q)_{\infty}}\sum\limits_{n\geq 0}(-1)^nq^{\frac{An(n+1)}{2}+a(n+1)}.
\end{align*}
\end{proof}
\subsection{Proof of Theorem \ref{mex3132}}
\begin{proof}
Using Theorem \ref{mexgen1}, we have
\begin{align*}
F_{3,1}(q)=\sum\limits_{n\geq 0}p_{3,1}(n)q^n=\frac{1}{(q)_{\infty}}\sum\limits_{n\geq 0}(-1)^nq^{\frac{n(3n-1)}{2}},
\\
F_{3,2}(q)=\sum\limits_{n\geq 0}p_{3,2}(n)q^n=\frac{1}{(q)_{\infty}}\sum\limits_{n\geq 0}(-1)^nq^{\frac{n(3n+1)}{2}}.
\end{align*}
Adding the two, we have
\begin{align*}
F_{3,1}(q)+F_{3,2}(q)&=\sum\limits_{n\geq 0}(p_{3,1}(n)+p_{3,2}(n))q^n
\\
&=\frac{1}{(q)_{\infty}}\sum\limits_{n\geq 0}(-1)^n(q^{\frac{n(3n-1)}{2}}+q^{\frac{n(3n+1)}{2}})
\\
&=\frac{1}{(q)_{\infty}}\Bigl(2+\sum\limits_{n\geq 1}(-1)^n(q^{\frac{n(3n-1)}{2}}+q^{\frac{n(3n+1)}{2}})\Bigr)
\\
&=1+\sum\limits_{n\geq 0}p(n)q^n,
\end{align*}
where the last step follows from Euler's pentagonal number theorem \cite[p.\ 11]{An98}.

Therefore, $\sum\limits_{n\geq 1}(p_{3,1}(n)+p_{3,2}(n))q^n=\sum\limits_{n\geq 1}p(n)q^n$. Comparing the coefficients of $q^n$ on both sides, we obtain the theorem.
\end{proof}
\subsection{Proof of Theorem \ref{mexrecur}}
\begin{proof}
Theorem \ref{mexgen1} gives 
\begin{align*}
F_{A,a}(q)&=\sum\limits_{n\geq 0}p_{A,a}(n)q^n
\\
&=\frac{1}{(q)_{\infty}}\sum\limits_{n\geq 0}(-1)^nq^{\frac{An(n-1)}{2}+an}
\\
&=\Bigl(\sum\limits_{n\geq 0}p(n)q^n\Bigr)\Bigl(1+\sum\limits_{n\geq 1}(q^{An(2n-1)+a.2n}-q^{A(n-1)(2n-1)+a(2n-1)})\Bigr)
\\
&=\sum\limits_{n\geq 0}\Bigl(p(n)+\sum\limits_{m\geq 1}(p(n-(A{2m \choose 2}+a(2m)))-
\\
&\qquad\qquad\qquad\qquad\quad p(n-(A{2m-1 \choose 2}+a(2m-1))))\Bigr)q^n.
\end{align*}
Comparing the coefficients of $q^n$ on both sides, we obtain the theorem.
\end{proof}
\begin{lemma}\label{PAap1}
If $n$ is a non-negative integer, then $p_{A,a}(n)=p(n)$ for every $a > n$.  
\end{lemma}
\begin{proof}
If $a>n$, then $a$ is not a part of every partition $\pi$ of $n$ and hence, we have mex$_{A,a}(\pi) = a$. The lemma follows. 
\end{proof}
\begin{corollary}\label{PAap2}
If $n$ is a non-negative integer, then $\overline{p}_{A,a}(n)=0$ for every $a > n$.
\end{corollary}
\subsection{Proof of Theorem \ref{mexrank1}}
\begin{proof}
Dyson \cite{Dy44} gave the generating function for the number of partitions of $n$ with rank $\geq j$.
\begin{align*}
\sum\limits_{m\geq j}\sum\limits_{n\geq 0}N(m,n)q^n&=\frac{1}{(q)_{\infty}}\sum\limits_{n\geq 1}(-1)^{n-1}q^{\frac{n(3n-1)}{2}}(1-q^n)\sum\limits_{m\geq j}q^{nm}
\\
&=\frac{1}{(q)_{\infty}}\sum\limits_{n\geq 1}(-1)^{n-1}q^{\frac{n(3n-1)}{2}+nj}
\\
&=\frac{1}{(q)_{\infty}}\sum\limits_{n\geq 0}(-1)^{n}q^{\frac{3n(n+1)}{2}+(n+1)(j+1)}
\\
&=\sum\limits_{n\geq 0}\overline{p}_{3,j+1}(n)q^n.
\end{align*}
Comparing the coefficients of $q^n$ on both sides, we obtain the theorem.
\end{proof}
As a result, we obtain the following corollaries.
\begin{corollary}\label{GardenofEden1}
If $n$ is a non-negative integer, then $\overline{p}_{3,3}(n)$ equals the number of Garden of Eden partitions of $n$.
\end{corollary}
\begin{proof}
The result follows from substituting $j=2$ in Theorem \ref{mexrank1}. Then, since $N(m,n)=N(-m,n)$, we have that $\overline{p}_{3,3}(n)$ equals the number of partitions of $n$ with rank $\geq 2$, which equals the number of partitions of $n$ with rank $\leq -2$, which is equinumerous with the Garden of Eden partitions. 
\end{proof}
\begin{corollary}\label{mexrank2}
If $n$ is a non-negative integer, then $p_{3,j+1}(n)$ equals the number of partitions of $n$ with rank $<j$, where j is a non-negative integer. 
\end{corollary}
\subsection{Proof of Theorem \ref{mexcrank1}}
\begin{proof}
Garvan \cite{Ga88} gave the generating function for the number of partitions of $n$ with crank $\geq j$.
\begin{align*}
\sum\limits_{m\geq j}\sum\limits_{n\geq 0}M(m,n)q^n&=\frac{1}{(q)_{\infty}}\sum\limits_{n\geq 1}(-1)^{n-1}q^{\frac{n(n-1)}{2}}(1-q^n)\sum\limits_{m\geq j}q^{nm}
\\
&=\frac{1}{(q)_{\infty}}\sum\limits_{n\geq 1}(-1)^{n-1}q^{\frac{n(n-1)}{2}+nj}
\\
&=\frac{1}{(q)_{\infty}}\sum\limits_{n\geq 0}(-1)^{n}q^{\frac{n(n+1)}{2}+(n+1)j}
\\
&=\sum\limits_{n\geq 0}\overline{p}_{1,j}(n)q^n.
\end{align*}
Comparing the coefficients of $q^n$ on both sides, we obtain the theorem.
\end{proof}
As a result, we obtain the following corollary.
\begin{corollary}\label{mexcrank2}
If $n$ is a non-negative integer, then $p_{1,j}(n)$ equals the number of partitions of $n$ with crank $<j$, where j is a non-negative integer. 
\end{corollary}
\subsection{Proof of Theorem \ref{AnNew1}}
\begin{proof}
The result follows from substituting $j=1$ in Corollary \ref{mexcrank2}. Then $p_{1,1}(n)$ equals the number of partitions of $n$ with crank $\leq 0$, which equals the number of partitions of $n$ with crank $\geq 0$ since $M(m,n)=M(-m,n)$.
\end{proof}
\subsection{Proof of Theorem \ref{AnNew2}}
\begin{proof}
The result follows from substituting $j=2$ in Corollary \ref{mexrank2}. Then $p_{3,3}(n)$ equals the number of partitions of $n$ with rank $\leq 1$, which equals the number of partitions of $n$ with rank $\geq -1$ since $N(m,n)=N(-m,n)$.
\end{proof}
\subsection{Proof of Theorem \ref{mexrankmom}}
\begin{proof}
We have the following generating function for the second rank moment due to Andrews \cite{An08} 
\begin{align*}
\sum\limits_{n\geq 0}\frac{1}{2} N_2(n)q^n&=\frac{-1}{(q)_{\infty}}\sum\limits_{n\geq 1}\dfrac{(-1)^nq^{\frac{n(3n+1)}{2}}(1+q^n)}{(1-q^n)^2}
\\
&=\frac{-1}{(q)_{\infty}}\sum\limits_{n\geq 1}(-1)^nq^{\frac{n(3n+1)}{2}}\sum\limits_{r\geq 0}(2r+1)q^{rn}
\\
&=\frac{-1}{(q)_{\infty}}\sum\limits_{r\geq 0}(2r+1)\sum\limits_{n\geq 1}(-1)^nq^{\frac{n(3n+1)}{2}+rn}
\\
&=\sum\limits_{r\geq 0}(2r+1)\left(\frac{1}{(q)_{\infty}}\sum\limits_{n\geq 0}(-1)^nq^{\frac{3n(n+1)}{2}+(r+2)(n+1)}\right)
\\
&=\sum\limits_{n\geq 0}\Bigl(\sum\limits_{r=0}^{n-2}(2r+1)\overline{p}_{3,r+2}(n)\Bigr)q^n,
\end{align*}
where the last step follows from Corollary \ref{PAap2}.

Similarly, the generating function for the second crank moment due to Garvan \cite{Ga11} is
\begin{align*}
\sum\limits_{n\geq 0}\frac{1}{2} M_2(n)q^n&=\frac{-1}{(q)_{\infty}}\sum\limits_{n\geq 1}\dfrac{(-1)^nq^{\frac{n(n+1)}{2}}(1+q^n)}{(1-q^n)^2}
\\
&=\frac{-1}{(q)_{\infty}}\sum\limits_{n\geq 1}(-1)^nq^{\frac{n(n+1)}{2}}\sum\limits_{r\geq 0}(2r+1)q^{rn}
\\
&=\frac{-1}{(q)_{\infty}}\sum\limits_{r\geq 0}(2r+1)\sum\limits_{n\geq 1}(-1)^nq^{\frac{n(n+1)}{2}+rn}
\\
&=\sum\limits_{r\geq 0}(2r+1)\left(\frac{1}{(q)_{\infty}}\sum\limits_{n\geq 0}(-1)^nq^{\frac{n(n+1)}{2}+(r+1)(n+1)}\right)
\\
&=\sum\limits_{n\geq 0}\Bigl(\sum\limits_{r=0}^{n-1}(2r+1)\overline{p}_{1,r+1}(n)\Bigr)q^n,
\end{align*}
where the last step follows from Corollary \ref{PAap2}. Comparing the coefficients of $q^n$ on both sides of the two equations, we obtain the theorem.
\end{proof}
\subsection{Proof of Theorem \ref{mexspt}}
\begin{proof}
Andrews \cite{An08} gave the following identity relating the smallest parts function to the second rank moment  $$\text{spt}(n)=np(n)-\frac{1}{2}N_2(n).$$
Dyson \cite{Dy89} showed that $\frac{1}{2}M_2(n)=np(n)$. Then, using Theorem \ref{mexrankmom}, we have 
\begin{align*}
\text{spt}(n)=np(n)-\frac{1}{2}N_2(n)&=\frac{1}{2}M_2(n)-\frac{1}{2}N_2(n)
\\
&=\sum\limits_{r=0}^{n-1}(2r+1)(\overline{p}_{1,r+1}(n)-\overline{p}_{3,r+2}(n))
\\
&=\sum\limits_{r=0}^{n-1}(2r+1)(p_{3,r+2}(n)-p_{1,r+1}(n)).
\end{align*}
\end{proof}
The rest of our results concern the connections of the extended mex function $p_{A,a}(n)$ with partitions of the integer $n$ having certain congruence conditions on the parts. As a motivation, we consider the celebrated Ramanujan theta function $$\psi(-q)=\sum\limits_{n=-\infty}^{\infty}(-1)^nq^{n(2n-1)}=\frac{(q^2;q^2)_{\infty}}{(-q;q^2)_{\infty}}.$$
Then \begin{align*}
(q^4;q^4)_{\infty}(q^3;q^4)_{\infty}(q;q^4)_{\infty}&=(q^4;q^4)_{\infty}(q;q^2)_{\infty}
\\
&=\frac{(q^2;q^2)_{\infty}(q;q^2)_{\infty}}{(q^2;q^4)_{\infty}}
\\
&=\frac{(q^2;q^2)_{\infty}}{(-q;q^2)_{\infty}}
\\
&=\sum\limits_{n=-\infty}^{\infty}(-1)^nq^{n(2n-1)}
\\
&=\sum\limits_{n\geq 0}(-1)^nq^{n(2n-1)}-\sum\limits_{n\geq 0}(-1)^nq^{2(n+1)^2+(n+1)}
\\
&=\sum\limits_{n\geq 0}(-1)^n(q^{\frac{4n(n-1)}{2}+n}-q^{\frac{4n(n+1)}{2}+3(n+1)})
\end{align*}
and so, $$\frac{1}{(q)_{\infty}}(q^4;q^4)_{\infty}(q^3;q^4)_{\infty}(q;q^4)_{\infty}=\sum\limits_{n\geq 0}(p_{4,1}(n)-\overline{p}_{4,3}(n))q^n.$$
Thus, we see that $p_{4,1}(n)-\overline{p}_{4,3}(n)$ equals the number of partitions of $n$ into parts congruent to $2$ modulo $4$. This find indicates the presence of a general class of identities involving congruences of similar type. We discover the two such families of identities in Theorem \ref{mexcongJTP} using Jacobi's triple product identity.
\begin{theorem}(Jacobi's triple product identity, \cite[p.\ 21]{An98})\label{JTP}
For $z\neq 0$ and $|q|<1$,
$$\sum\limits_{n=-\infty}^{\infty}z^nq^{n^2}=\prod\limits_{n\geq 0}(1-q^{2n+2})(1+zq^{2n+1})(1+z^{-1}q^{2n+1}).$$
\end{theorem}
\begin{theorem}(Andrews, \cite[p.\ 22]{An98})\label{JTP1}
The following is an equivalent form of Jacobi's triple product identity, due to Andrews
$$\sum\limits_{n=-\infty}^{\infty}(-1)^nq^{(2k+1)n(n+1)/2-in}=\prod\limits_{n\geq 0}(1-q^{(2k+1)(n+1)})(1-q^{(2k+1)n+i})(1-q^{(2k+1)(n+1)-i}).$$
\end{theorem}
\begin{lemma}\label{JTP2}
The following is another equivalent form of Jacobi's triple product identity:
$$\sum\limits_{n=-\infty}^{\infty}(-1)^nq^{kn(n+1)-in}=\prod\limits_{n\geq 0}(1-q^{2k(n+1)})(1-q^{2kn+i})(1-q^{2k(n+1)-i}).$$
\end{lemma}
\begin{proof}
The proof follows from substituting $z \rightarrow -q^{k-i}$ and $q \rightarrow q^k$ in Jacobi's triple product identity in Theorem \ref{JTP}.
\end{proof}
\subsection{Proof of Theorem \ref{mexcongJTP}}
\begin{proof}
\begin{enumerate}[(a)]
\item Using Lemma \ref{JTP2}, we have 
\begin{align*}
\prod\limits_{\substack{n\geq 1 \\ n \not\equiv 0, \pm i\text{ (mod }2k)}}\frac{1}{1-q^n}&=\frac{1}{(q)_{\infty}}\sum\limits_{n\geq 0}(-1)^nq^{\frac{2kn(n+1)}{2}-in}(1-q^{(2n+1)i})
\\
&=\frac{1}{(q)_{\infty}}\sum\limits_{n\geq 0}(-1)^nq^{\frac{2kn(n-1)}{2}+(2k-i)n}-\frac{1}{(q)_{\infty}}\sum\limits_{n\geq 0}(-1)^nq^{\frac{2kn(n+1)}{2}+i(n+1)}
\\
&=\sum\limits_{n\geq 0}(p_{2k,2k-i}(n)-\overline{p}_{2k,i}(n))q^n.
\end{align*}
\item Using Theorem \ref{JTP1}, we have 
\begin{align*}
\prod\limits_{\substack{n\geq 1 \\ n \not\equiv 0, \pm i\text{ (mod }2k+1)}}\frac{1}{1-q^n}&=\frac{1}{(q)_{\infty}}\sum\limits_{n\geq 0}(-1)^nq^{\frac{(2k+1)n(n+1)}{2}-in}(1-q^{(2n+1)i})
\\
&=\frac{1}{(q)_{\infty}}\sum\limits_{n\geq 0}(-1)^n(q^{\frac{(2k+1)n(n-1)}{2}+(2k+1-i)n}-q^{\frac{(2k+1)n(n+1)}{2}+i(n+1)})
\\
&=\sum\limits_{n\geq 0}(p_{2k+1,2k+1-i}(n)-\overline{p}_{2k+1,i}(n))q^n.
\end{align*}
\end{enumerate}
Comparing the coefficients of $q^n$ on both sides in (a) and (b) above, we obtain the theorem.
\end{proof}
\subsection{Proof of Theorem \ref{AnNew3}}
\begin{proof}
The following is a consequence of Cauchy's identity due to Euler \cite[p.\ 19]{An98} 
$$\sum\limits_{n\geq 0}\frac{t^n}{(q)_{n}}=\frac{1}{(t;q)_{\infty}}.$$
Substituting $t=-q$ in the above identity, we get 
\begin{align*}
&\quad\sum\limits_{n\geq 0}\frac{(-1)^nq^n}{(q)_{n}}=\frac{1}{(-q;q)_{\infty}},
\\
&\quad\sum\limits_{n\geq 0}\frac{q^{2n}}{(q)_{2n}}-\sum\limits_{n\geq 0}\frac{q^{2n+1}}{(q)_{2n+1}}=(q;q^2)_{\infty},
\end{align*}
\begin{equation}\label{peo}
\sum\limits_{n\geq 0}(p_e(n)-p_o(n))q^n=(q;q^2)_{\infty}.
\end{equation}
Theorem \ref{mexcongJTP} (a) gives
\begin{align*}
\prod\limits_{\substack{n\geq 1 \\ n \not\equiv 0, \pm i\text{ (mod }2k)}}\frac{1}{1-q^n}&=\sum\limits_{n\geq 0}(p_{2k,2k-i}(n)-\overline{p}_{2k,i}(n))q^n.
\end{align*}
Multiplying by $(q)_{\infty}$ on both sides, we have
\begin{align*}
\prod\limits_{\substack{n\geq 0}}(1-q^{2k(n+1)})(1-q^{2kn+i})(1-q^{2k(n+1)-i})&=(q)_{\infty}\sum\limits_{n\geq 0}(p_{2k,2k-i}(n)-\overline{p}_{2k,i}(n))q^n.
\end{align*}
Using Lemma \ref{JTP2}, we then have
\begin{align*}
\sum\limits_{n=-\infty}^{\infty}(-1)^nq^{kn(n+1)-in}=(q)_{\infty}\sum\limits_{n\geq 0}(p_{2k,2k-i}(n)-\overline{p}_{2k,i}(n))q^n.
\end{align*}
Substituting $k=i=1$, we get
\begin{align*}
\sum\limits_{n=-\infty}^{\infty}(-1)^nq^{n^2}=(q)_{\infty}\sum\limits_{n\geq 0}(p_{2,1}(n)-\overline{p}_{2,1}(n))q^n,
\end{align*}
Then, using Gauss' identity \cite[Eq.\ (2.2.12), p.\ 23]{An98}, we have
\begin{equation}\label{p21}
\sum\limits_{n\geq 0}(p_{2,1}(n)-\overline{p}_{2,1}(n))q^n=\frac{1}{(q)_{\infty}}\sum\limits_{n=-\infty}^{\infty}(-1)^nq^{n^2}=\dfrac{1}{(-q)_{\infty}}=(q;q^2)_{\infty}.
\end{equation}
Comparing equations \eqref{peo} and \eqref{p21}, we have $p_e(n)-p_o(n) = p_{2,1}(n)-\overline{p}_{2,1}(n)$. Also, $p_e(n)+p_o(n) = p(n) = p_{2,1}(n)+\overline{p}_{2,1}(n)$. Adding and subtracting, we have $p_{2,1}(n)=p_e(n)$ and $\overline{p}_{2,1}(n)=p_o(n)$ respectively.
\end{proof}
The three results below allow us to obtain shifted partition type identities between the difference of mexes and partitions with congruence conditions on the parts. 
\begin{theorem}(Blecksmith, Brillhart, and Gerst \cite{Bl-Br-Ge87})\label{congid1}
\begin{enumerate}[(a)]
\item $$\prod\limits_{\substack{n\geq 1 \\ n \not\equiv \pm (2,8,12,14){\rm \,\,(mod}\,\, 32{\rm )}}} (1-q^n)=\sum\limits_{n\geq 1}(-1)^n(q^{2n^2-1}-q^{n^2-1}).$$
\item $$\prod\limits_{\substack{n\geq 1 \\ n \not\equiv \pm (1,4,6,8,10,11){\rm \,\,(mod}\,\, 24{\rm )}}} (1-q^n)=\sum\limits_{n\geq 1}(-1)^n(q^{3n^2-1}-q^{n^2-1}).$$
\end{enumerate}
\end{theorem}
\begin{theorem}(Blecksmith, Brillhart, and Gerst \cite{Bl-Br-Ge88})\label{congid2}
$$\prod\limits_{\substack{n\geq 1 \\ n \equiv 0,\pm 3 \\ {\rm (mod}\,\, 10{\rm )}}} (1-q^n)
\prod\limits_{\substack{n\geq 1 \\ n \equiv \pm 4 \\ {\rm (mod}\,\, 40{\rm )}}} (1-q^n)
\prod\limits_{\substack{n\geq 1 \\ n \equiv \pm 8 \\ {\rm (mod}\,\, 20{\rm )}}} (1+q^n)=\sum\limits_{n\geq 1}(-1)^n(q^{5n^2-1}-q^{n^2-1}).$$
\end{theorem}
\subsection{Proof of Theorem \ref{mexcong1}}
\begin{proof}
\begin{align*}
\sum\limits_{n\geq 1}(-1)^n(q^{2n^2-1}-q^{n^2-1})&=\sum\limits_{n\geq 0}(-1)^nq^{n^2+2n}-\sum\limits_{n\geq 0}(-1)^nq^{2n^2+4n+1}
\\
&=(q)_{\infty}\sum\limits_{n\geq 0}p_{2,3}(n)q^n-(q)_{\infty}\sum\limits_{n\geq 0}p_{4,6}(n)q^{n+1}
\\
&=(q)_{\infty}\Bigl(1+\sum\limits_{n\geq 1}(p_{2,3}(n)-p_{4,6}(n-1))q^n\Bigr)
\\
&=(q)_{\infty}\sum\limits_{n\geq 0}(p_{2,3}(n)-p_{4,6}(n-1))q^n,
\end{align*}
where the last step follows from the fact that $p_{A,a}(0)=1$ and $p_{A,a}(n)=0$ for $n<0$. Then, dividing both sides by $(q)_{\infty}$ and using Theorem \ref{congid1} (a), we have our result.
\end{proof}
\subsection{Proof of Theorem \ref{mexcong2}}
\begin{proof}
\begin{align*}
\sum\limits_{n\geq 1}(-1)^n(q^{3n^2-1}-q^{n^2-1})&=\sum\limits_{n\geq 0}(-1)^nq^{n^2+2n}-\sum\limits_{n\geq 0}(-1)^nq^{3n^2+6n+2}
\\
&=(q)_{\infty}\sum\limits_{n\geq 0}p_{2,3}(n)q^n-(q)_{\infty}\sum\limits_{n\geq 0}p_{6,9}(n)q^{n+2}
\\
&=(q)_{\infty}\Bigl(1+\sum\limits_{n\geq 1}(p_{2,3}(n)-p_{6,9}(n-2))q^n\Bigr)
\\
&=(q)_{\infty}\sum\limits_{n\geq 0}(p_{2,3}(n)-p_{6,9}(n-2))q^n,
\end{align*}
where the last step follows from the fact that $p_{A,a}(0)=1$ and $p_{A,a}(n)=0$ for $n<0$. Then, dividing both sides by $(q)_{\infty}$ and using Theorem \ref{congid1} (b), we have our result.
\end{proof}
\subsection{Proof of Theorem \ref{mexcong3}}
\begin{proof}
\begin{align*}
\sum\limits_{n\geq 1}(-1)^n(q^{5n^2-1}-q^{n^2-1})&=\sum\limits_{n\geq 0}(-1)^nq^{n^2+2n}-\sum\limits_{n\geq 0}(-1)^nq^{5n^2+10n+4}
\\
&=(q)_{\infty}\sum\limits_{n\geq 0}p_{2,3}(n)q^n-(q)_{\infty}\sum\limits_{n\geq 0}p_{10,15}(n)q^{n+4}
\\
&=(q)_{\infty}\Bigl(1+\sum\limits_{n\geq 1}(p_{2,3}(n)-p_{10,15}(n-4))q^n\Bigr)
\\
&=(q)_{\infty}\sum\limits_{n\geq 0}(p_{2,3}(n)-p_{10,15}(n-4))q^n,
\end{align*}
where the last step follows from the fact that $p_{A,a}(0)=1$ and $p_{A,a}(n)=0$ for $n<0$. Then, dividing both sides by $(q)_{\infty}$ and using Theorem \ref{congid2} where we convert all the congruence conditions to modulo $40$, we have
\begin{align*}
\sum\limits_{n\geq 0}(p_{2,3}(n)-p_{10,15}(n-4))q^n&=\dfrac{\prod\limits_{\substack{n\geq 1 \\ n \equiv 0,\pm (3,4,7,10,13,17),20 \\ {\rm (mod}\,\, 40{\rm )}}} (1-q^n)\prod\limits_{\substack{n\geq 1 \\ n \equiv \pm (8,12) \\ {\rm (mod}\,\, 40{\rm )}}} (1+q^n)}{(q)_{\infty}}
\\
&=\dfrac{\prod\limits_{\substack{n\geq 1 \\ n \equiv \pm (8,12) \\ {\rm (mod}\,\, 40{\rm )}}} (1+q^n)}{\prod\limits_{\substack{n\geq 1 \\ n \not\equiv 0,\pm (3,4,7,10,13,17),20 \\ {\rm (mod}\,\, 40{\rm )}}} (1-q^n)}.
\end{align*}
Comparing the coefficients of $q^n$ on both sides, we obtain the theorem.
\end{proof}
\section{Some auxiliary results}\label{auxiliary}
We present some interesting properties of the extended mex function $p_{A,a}(n)$ that we came across through computational evidence. 
\begin{theorem}\label{ppbar}
If $n$ is a non-negative integer, then $p_{A,a}(n-(a-A))=\overline{p}_{A,a-A}(n)$ for every $a > A$.
\end{theorem}
\begin{proof}
For all $a > A$, we have 
\begin{align*}
\sum\limits_{n\geq 0}\overline{p}_{A,a-A}(n)q^n&=\frac{1}{(q)_{\infty}}\sum\limits_{n\geq 0}(-1)^nq^{\frac{An(n+1)}{2}+(a-A)(n+1)}
\\
&=\frac{1}{(q)_{\infty}}\sum\limits_{n\geq 0}(-1)^nq^{\frac{An(n-1)}{2}+an+(a-A)}
\\
&=\sum\limits_{n\geq 0}{p}_{A,a}(n)q^{n+(a-A)}
\\
&=\sum\limits_{n\geq 0}{p}_{A,a}(n-(a-A))q^n.
\end{align*}
Comparing the coefficients of $q^n$ on both sides, we obtain the theorem.
\end{proof}
\begin{corollary}\label{GardenofEden2}
If $n$ is a non-negative integer, then ${p}_{3,6}(n-3)$ equals the number of Garden of Eden partitions of $n$.
\end{corollary}
\begin{proof}
The result follows from Corollary \ref{GardenofEden1} and considering $A=3$, $a=6$ in Theorem \ref{ppbar}.
\end{proof}
\begin{theorem}\label{pnnn1}
For any integer $n \geq 2,\, p_{n,k}(n)=p(n)-p(n-k)$ for every integer $k\ge 1$.
\end{theorem}
\begin{proof}
$p_{n,k}(n)$ enumerates all partitions of $n$ except $n-k$ of those whose largest part is less than or equal to $k$. Hence $p_{n,k}(n) = p(n)-p(n-k)$.
\end{proof}
\begin{theorem}\label{p31}
For any non-negative integer $n$, we have
\begin{align*}
p_{3,n}(n)-p_{1,n-1}(n)=0\, \forall\, n \geq 2,
\\
p_{3,n+1}(n)-p_{1,n}(n)=1\, \forall\, n \geq 1.
\end{align*}
\end{theorem}
\begin{proof}
The statements can be proved using the recurrence relation in Theorem \ref{mexrecur} similar to the proof of Lemma \ref{PAap1}.
\end{proof}
\section{Minimal excludant of overpartitions}\label{overmex}
Analogous to Andrews' and Newman's minimal excludant or ``mex" of a partition function, we define what we call the ``mex" of overpartitions. A surprising but welcome discovery is its connection to the Ramanujan function $R(q)$. This is a work in progress and will appear in a forthcoming paper. We present our initial findings here.

Recall that an overpartition is a partition in which the first occurrence
of a number may be overlined \cite{Cor-Love04}. For an overpartition $\pi$, we define $\overline{\text{mex}}(\pi)$ to be the smallest positive integer that is not a part of $\pi$. For example, if we consider the following two overpartitions
\begin{align*}
\pi_1&=\overline{14}+12+9+7+\overline{4}+3+\overline{2}+1,
\\
\pi_2&=\overline{5}+\overline{4}+4+2+1,
\end{align*}
then, $\overline{\text{mex}}(\pi_1)=5$ and $\overline{\text{mex}}(\pi_2)=3$. Let $\overline{m}(n)$ denote the number of overpartitions of $n$ having the property that no positive integer less than  $\overline{\text{mex}}(\pi)$ is overlined. As an example, below are the eight overpartitions of $3$
\begin{align*}
3, \overline{3}, 2+1, \overline{2}+1, 2+\overline{1}, \overline{2}+\overline{1}, 1+1+1, \overline{1}+1+1.
\end{align*}

Then, for $\pi \in \{3, \overline{3}, 2+1, 1+1+1\},$ the
quantity $\overline{\text{mex}}(\pi)$ has the stated property, and hence $\overline{m}(3)=4$.
Let $\overline{M}(q)=\sum\limits_{n=0}^{\infty}\overline{m}(n)q^n$. Then
\begin{theorem}\label{overptnmex}
We have $\overline{M}(q)=\overline{P}(q)\,(2-R(q))$, where $\overline{P}(q)=\frac{(-q;q)_{\infty}}{(q;q)_{\infty}}$ is the generating function for the number of overpartitions of $n$ and $R(q)=\sum\limits_{n=0}^{\infty}\dfrac{q^{\frac{n(n+1)}{2}}}{(-q;q)_n}$ is Ramanujan's function from his Lost Notebook.
\end{theorem}
\begin{remark}
The function $R(q)$ has been studied by Andrews \cite{An86}. It is also the generating function for number of partitions into distinct parts with even rank minus those with odd rank.
\end{remark}
To prove the above theorem, we need a result of Gasper and Rahman \cite[Appendix II, Formula II.5]{Ga-Rah04}, which we state as a lemma below.
\begin{lemma}
${}_1 \phi_1(a;c;q,\frac{c}{a})=\dfrac{(\frac{c}{a};q)_{\infty}}{(c;q)_{\infty}}$.
\end{lemma} 
The proof of the lemma follows easily from a suitable substitution and limit in Heine's transformation. We replace $a\mapsto q$ and $c\mapsto -zq$ in the above lemma to get 
\begin{equation}\label{gasrah}
\sum\limits_{n=0}^{\infty}\frac{z^nq^{\frac{n(n-1)}{2}}}{(-zq;q)_n}=1+z.
\end{equation}
\subsection{Proof of Theorem \ref{overptnmex}}
\begin{proof}
By standard combinatorial arguments, we deduce that
\begin{align*}
\overline{M}(q)&=\sum\limits_{n=0}^{\infty}\overline{m}(n)q^n
\\
&=\sum\limits_{n=0}^{\infty}\dfrac{q^{1+2+\cdots+(n-1)}\prod\limits_{m=n+1}^{\infty}(1+q^m)}{\prod\limits_{\substack{m=1 \\ m \neq n}}^{\infty}(1-q^m)}
\\
&=\frac{(-q;q)_{\infty}}{(q;q)_{\infty}}\sum\limits_{n=0}^{\infty}\dfrac{q^{\frac{n(n-1)}{2}}(1-q^n)}{(-q;q)_n}
\\
&=\frac{(-q;q)_{\infty}}{(q;q)_{\infty}}\left(\sum\limits_{n=0}^{\infty}\dfrac{q^{\frac{n(n-1)}{2}}}{(-q;q)_n}-\sum\limits_{n=0}^{\infty}\dfrac{q^{\frac{n(n+1)}{2}}}{(-q;q)_n}\right)
\\
&=\overline{P}(q)(2-R(q))\hspace{3mm}(\text{substituting}~z=1~\text{in}~\text{equation}~\eqref{gasrah}).
\end{align*}
\end{proof}
\section{Conclusion}
It would be interesting to investigate more identities of the shifted partition nature as in Theorems \ref{mexcong1}, \ref{mexcong2}, and \ref{mexcong3}, since a pattern emerges in the difference of the extended mex function in the three identities. A related line of inquiry would be if certain other linear combination of mexes lead to further restricted partition identities. Secondly, the role and significance of the coefficient $2r+1$ in the identity for \text{spt}($n$) sparks our curiosity. Thirdly, keeping in tradition with a very common question raised in the literature on the minimal excludant, bijective proofs of the theorems stated in the paper would be of great interest. Lastly, although our major focus in the paper has been connecting the mex with other partition statistics, the last two auxiliary results indicate the existence of more such properties within the extended mex function $p_{A,a}(n)$ and could be a realm worthy of further exploration.
\section{Acknowledgments}
We express our gratitude to George Andrews and Frank Garvan for previewing a preliminary draft of this paper and their very helpful comments and suggestions.

\bigskip
\hrule
\bigskip

\noindent 2020 {\it Mathematics Subject Classification}: Primary 11P81; Secondary 11P82, 11P83, 05A17.

\noindent \emph{Keywords:} minimal excludant, extended mex function, partition statistic, rank, crank, spt, partition congruence, overpartition.

\bigskip
\hrule
\bigskip

\vspace*{+.1in}
\noindent
Received October 21 2022;
revised version received February 9 2023; March 16 2023. 
Published in {\it Journal of Integer Sequences}, March 17 2023.

\bigskip
\hrule
\bigskip

\noindent
Return to
\href{https://cs.uwaterloo.ca/journals/JIS/}{Journal of Integer Sequences home page}.
\vskip .1in


\begin{thebibliography}{99}
\bibitem{An86}
G.~E. Andrews, Ramanujan's Lost Notebook V:Euler's partition identity, \emph{Adv. Math.} {\bf 61} (1986),
156--164.
\bibitem{An98}
G. E. Andrews, \emph{The Theory of Partitions}, Cambridge University Press, 1998.
\bibitem{An08}
G. E. Andrews, The number of smallest parts in the partitions of $n$, \emph{J. Reine Angew. Math.}
{\bf 624} (2008),  133--142.
\bibitem{An-New19}
G. E. Andrews and D.~Newman, Partitions and the minimal excludant, \emph{Ann. Comb.}
{\bf 23} (2019), 249--254.
\bibitem{An-New20}
G. E. Andrews and D. Newman, The minimal excludant in integer partitions, \emph{J. Integer Sequences}
{\bf 23} (2020), \href{https://cs.uwaterloo.ca/journals/JIS/VOL23/Andrews/andrews5.html}{Article 20.2.3}.
\bibitem{Atk-Ga03}
A. O. L. Atkin and F. G. Garvan, Relations between the ranks and cranks of partitions, \emph{Ramanujan J.}
{\bf 7} (2003), 343--366.
\bibitem{Bl-Br-Ge87}
R. Blecksmith, J. Brillhart, and I. Gerst, Parity results for certain partition functions and identities similar to theta function identities, \emph{Math. Comp.} {\bf 48} (1987), 29--38.
\bibitem{Bl-Br-Ge88}
R. Blecksmith, J. Brillhart, and I. Gerst, Some infinite product identities, \emph{Math. Comp.} {\bf 51} (1988), 301--314.
\bibitem{Cor-Love04}
S. Corteel and J. Lovejoy, Overpartitions, \emph{Trans. Amer. Math. Soc.} {\bf 356} (2004), 1623--1635.
\bibitem{Dy44}
F. J. Dyson, Some guesses in the theory of partitions, \emph{Eureka} {\bf 8} (1944), 10--15.
\bibitem{Dy89}
F. J. Dyson, Mappings and symmetries of partitions, \emph{J. Combin. Theory Ser. A} {\bf 51} (1989),
169--180.
\bibitem{Ga88}
F. G. Garvan, New combinatorial interpretations of Ramanujan's partition congruences
mod $5$, $7$ and {$11$}, \emph{Trans. Amer. Math. Soc.} {\bf 305} (1988), 47--77.
\bibitem{Ga11}
F. G. Garvan, Higher order \text{spt}-functions, \emph{Adv. Math.} {\bf 228} (2011), 241--265.
\bibitem{Ga-Rah04}
G. Gasper and M. Rahman, \emph{Basic Hypergeometric Series}, 2nd ed., Cambridge University Press, 2004.
\bibitem{Hop-Sel07}
B. Hopkins and J. A. Sellers, Exact enumeration of {G}arden of {E}den partitions.
\newblock In B. Landman, M. Nathanson, J. Ne\v{s}et\v{r}il, R. Nowakowski, and C. Pomerance, eds., \emph{Combinatorial Number Theory}, 
De Gruyter, 2007, pp.\ 299--303. 
\end{thebibliography}
\end{document}